\documentclass[12pt]{article}
\usepackage{graphics}
\usepackage{amsmath}
\usepackage{amssymb}
\usepackage{hyperref}
\usepackage{graphicx}
\usepackage{makeidx}
\usepackage{mathrsfs}
\usepackage{amsfonts}
\usepackage[mathscr]{eucal}
\usepackage{amsmath}

\DeclareGraphicsExtensions{.eps,.pdf,.jpg,.png}

\def\RR{\vbox {\hbox to 8.9pt {I\hskip-2.1pt R\hfil}}}

\def\pni{\par\noindent}
\def\vsh{\smallskip}

\def\vsp{\vsh\pni} 

\numberwithin{equation}{section}

\begin{document}
\font\title=cmbx12 scaled\magstep2
\font\bfs=cmbx12 scaled\magstep1
\font\little=cmr10
\begin{center}
{\title Pulse waves in the viscoelastic  Kelvin-Voigt model:
a revisited  approach}
\\ [0.50 truecm]
Juan Luis GONZ{\'A}LEZ  SANTANDER$^{(1)}$
 \\  [0.25 truecm]
 Francesco MAINARDI$^{(2)}$, and
 Andrea MENTRELLI$^{(3) (4)}$
 \\[0.50 truecm]
 $^{(1)}$  Department of Mathematics, University of  Oviedo.
 \\ C Leopolodo Calvo Sotelo 18, 33007 Oviedo, Spain;
 {gonzalezmarjuan@uniovi.es}
 \\
$^{(2)}$
 Department of Physics and Astronomy, University of Bologna and INFN.
\\ {Via Irnerio 46, I-40126 Bologna, Italy};
\\{ francesco.mainardi@unibo.it; mainardi@bo.infn.it; fracalmo@gmail.com}
\\
$^{(3)}$ Department of Mathematics and AM$^2$,  University of Bologna.
\\ Via Saragozza 8, 40123 Bologna, Italy;
 {andrea.mentrelli@unibo.it}
\\
$^{(4)}$ Istituto Nazionale di Fisica Nucleare (I.N.F.N.), Sezione di Bologna,
\\  I.S. FLAG,  Viale Berti Pichat 6/2, 40127 Bologna, Italy.
\

 \vskip 0.25truecm
 {\bf Published in Mathematics (MDPI)
  Vol.14 No 3 (2025), 528/1--16}
  \\ {\bf DOI:
10.3390/math14030528}
\end{center}

\begin{abstract}
We calculate the mechanical response $r\left(x,t\right)$ of an initially quiescent semi-infinite homogeneous medium to a pulse applied at the origin, and this is achieved
within the framework of the Kelvin--Voigt model. Although this problem has been extensively studied in the literature because of its wide range of applications---particularly in seismology---here, we present a solution in a novel integral form. This integral solution avoids the numerical computation of the solution
in terms of the inverse Laplace transform; that is, numerical integration in the complex plane. In particular, we derive integral form expressions
for both delta-pulse and step-pulse excitations which are simpler and more computationally efficient than those previously reported in the literature.
Furthermore, the obtained expressions allow us to obtain simple asymptotic formulas for $r\left(x,t\right)$ as $x,t \rightarrow 0,\infty$ for both step- and {delta-type pulses}.

 \vsp {\bf Keywords}:
{Transient waves in linear viscoelasticity; Kelvin-Voigt model;
Laplace transform}

\vsp {\bf Mathematics Subject Classification (MSC)}:
{41A60,  44A10, 35L20, 33C10}

\end{abstract}
\newpage

%



\section{Introduction}

The problem of determining the response of uniaxial waves in an initially quiescent semi-infinite medium is fundamental in linear viscoelasticity and is discussed in most standard textbooks on rheology.
Indeed, it constitutes a preliminary step for the analysis of wave propagation in linear dispersive media subject to dissipation, including seismic waves in the Earth.
 Following the notation used in Mainardi's book \cite{MainardiBook},
 we consider the mechanical response $r\left( x,t\right) $ at time $t\geq 0 $ of an initially quiescent, semi-infinite homogeneous medium extending over $x\geq 0$ (with density $\rho$) to a pulse applied at the origin; that is, $r_{0}\left( t\right) =r\left( 0,t\right) $.
  As pointed out in Hunter's review paper \cite{Hunter}, the quantity $r(x,t)$ may represent the stress $\sigma(x,t)$, strain $\varepsilon(x,t)$, displacement $u(x,t)$, or particle velocity $v(x,t)$.
  According to the classical theory described in \cite{MainardiBook} (Sect. 4.2.1) (see also \cite{Hunter}),
  the response is obtained through the inversion of Laplace transform of the form 

  \begin{equation}
  r\left( x,t\right) =\frac{1}{2\pi i}\int_{\gamma-i \,\infty}^{\gamma+i \,\infty}\tilde{r}_{0}\left( s\right)
\exp \left( s \, t- m(s) \,x \right) \,ds,
\label{r(x,t)_inv-m_Laplace}
\end{equation}
  where $\gamma>0$, $s$ is the Laplace parameter, $\widetilde{r}_{0}\left( s\right) $ is the Laplace transform of the applied pulse $r_{0}\left( t\right) $,
and $m\left( s\right) $ is a characteristic complex function depending on the
Laplace transform of the material functions; that is, the creep compliance $J(t)$ (the strain to a step input of stress) and the relaxation modulus $G(t)$ (the stress to a step input of strain), according
to the equation:
\begin{equation}
m(s) = s \sqrt{\rho\, s\,\widetilde J(s)}  =
s \sqrt{\frac{\rho}{s\,\widetilde G(s)}}.
\label{m(s)JG}
\end{equation}
Here, we have used the notation given in \cite{MainardiBook} where $\widetilde J(s)$ and $\widetilde G(s)$
 represent the Laplace transforms of the creep compliance and relaxation modulus, respectively.

 At this point, we would distinguish the case of the existence---or lack thereof---of a finite wave front velocity $c>0$ according to
 \begin{equation}
 \lim_{|s| \to \infty} \frac{m(s)}{s} =
 \left \{
 \begin{array} {ll}
  0 \quad {\rm if}  \quad 1/{J(0)}= G(0)=\infty, \\
  \frac{1}{c} \quad  {\rm if} \quad 1/{J(0)}= G(0) = \rho\, c^2,
  \end{array}
  \right .
  \end{equation}
  so that  only in the second case we have a finite wave-front
   velocity
   \begin{equation}
   c= \sqrt{\frac{1}{J(0)\,\rho}} = \sqrt{\frac{G(0)}{\rho}}\,, \label{c_def}
   \end{equation}
and, consequently, we introduce the index of refraction $n(s)$ defined by
\begin{equation}
n(s) = c \, \frac {m(s)}{s}.
\label{index of refraction}
\end{equation}
In 
 this case
\begin{equation}
r\left( x,t\right) =\frac{1}{2\pi i}\int_{\gamma-i \,\infty}^{\gamma+i \,\infty}\tilde{r}_{0}\left( s\right)
\exp \left( s\left[ t-n\left( s\right) \,x/c\right] \right) \,ds,
\label{r(x,t)_inv_Laplace}
\end{equation}%
so that the response turns out to be  zero if $x> ct$, namely for
$t<x/c$. This justifies the definition of wave-front velocity.
If we do not have a finite wave front velocity, we have a diffusion phenomenon;
that is, an instantaneous propagation of the signal like in \mbox{heat conduction.}

In \cite{MainardiBook}, we found that the viscoelastic bodies are classified into four types
according to their instantaneous and equilibrium responses.
In fact, we easily recognize four
possibilities  for the limiting values of the creep
compliance and relaxation modulus, as listed in Table 1,
where we have denoted $J_g =J(0)$, $J_e=J(\infty)$,
$G_g=G(0)$, $G_e= G(\infty)$ (with $J_g\, G_g =J_e \, G_e =1$).
\vskip 0.25truecm
\begin{center}
\begin{tabular}{|c||c|c||c|c|}
\hline
Type & $J_g$ & $J_e$ & $G_g$ & $G_e$ \\
\hline
I   & $>0$ & $<\infty$ & $<\infty$ & $>0$ \\
II  & $>0$ & $=\infty$ & $<\infty$ & $=0$ \\
III & $=0$ & $<\infty$ & $=\infty$ & $>0$ \\
IV  & $=0$ & $=\infty$ & $=\infty$ & $=0$ \\
\hline
\end{tabular}  
\vskip 0.25truecm
{\bf Table 1} \ {The four types of viscoelasticity}.
\end{center}
\vsp
We note that the viscoelastic bodies of type I
exhibit   both instantaneous and equilibrium elasticity,
so their behavior appears  close to the purely elastic one
for sufficiently short and  long times.
The  bodies of  type II and IV exhibit
a complete stress relaxation (at constant strain) since $G_e =0$ and an
infinite strain creep (at constant stress) since $J_e = \infty\,,$
so  they do not show equilibrium elasticity.
However, the bodies of type III and IV
do not show instantaneous elasticity  since $J_g = 0\,$
($G_g =\infty$).
Concerning wave propagation, according to (\ref{c_def}), only bodies of type I and II exhibit a finite wave-front velocity $c$, whereas bodies of type III and IV exhibit diffusion; that is, $c \rightarrow\infty$.
The simplest body of type III is the Kelvin--Voigt model, which is the subject of this work.
The Kelvin--Voigt model is relevant to a wide range of transient wave problems,
 particularly in seismology, as discussed in numerous classical and modern studies
 \cite{Jeffreys,Lee,Morrison,Hanin,Collins,Clark,Clark2,%
 Jaramillo,Mainardi1972,Buchen1975,Mainardi1975,Dozio,Colombaro}.

 In this paper, we revisit the transient-wave problem in a Kelvin--Voigt medium and derive integral-form solutions for both step-pulse and delta-pulse excitations. We also analyze the asymptotic behavior of these solutions as $t\rightarrow 0,\infty$ and $x\rightarrow 0,\infty$. Our aim is to improve and extend upon previous results available in the literature.

\section{The Kelvin--Voigt Model}

\subsection{General Solution}

The constitutive relation in the Kelvin--Voigt model is%
\begin{equation}
\sigma =G_{e}\left[ \varepsilon +t_{\varepsilon }\frac{\partial \varepsilon
}{\partial t}\right] ,  \label{Eq_constitutive_Kelvin-Voigt}
\end{equation}%
where $t_{\varepsilon }>0$ is the retardation time and $G_{e}$ the
equilibrium\ modulus. Apply the Laplace transform to the
constitutive relation (\ref{Eq_constitutive_Kelvin-Voigt})\ to obtain%
\begin{equation}
\tilde{\sigma}\left( x,s\right) =G_{e}\,\left[ \,\tilde{\varepsilon}\left(
x,s\right) +t_{\varepsilon }\left( s\,\tilde{\varepsilon}\left( x,s\right)
-\varepsilon \left( x,0\right) \right) \right] .
\label{Eq_Kelvin-Voigt_Laplace_cond}
\end{equation}
If
 we consider that, initially, the rod is unstrained, we have
\begin{equation*}
\varepsilon \left( x,0\right) =0,
\end{equation*}%
thus, (\ref{Eq_Kelvin-Voigt_Laplace_cond})\ is reduced to%
\begin{equation}
\tilde{\sigma}\left( x,s\right) =G_{e}\left( 1\,+t_{\varepsilon }s\right) \,%
\tilde{\varepsilon}\left( x,s\right) .
\label{Eq_constitutive_Kelvin-Voigt_Laplace}
\end{equation}%
Apply to (\ref{Eq_constitutive_Kelvin-Voigt_Laplace})\ the derivative with
respect to $x$ and take into account the basic equations of linear viscoelasticity in the Laplace domain (see \cite{MainardiBook}, Sect. 4.2.1),
 i.e., the equation of motion:
\begin{equation}
\frac{\partial }{\partial x}\tilde{\sigma}\left( x,s\right) = \rho \,s^{2}\,\tilde{r}\left( x,s\right),  \label{Eq_motion_Laplace}
\end{equation}
and the kinematic equation:
\begin{equation}
\tilde{\varepsilon}\left( x,s\right) =\frac{\partial }{\partial x}\tilde{r}%
\left( x,s\right),  \label{Eq_kinematic_Laplace}
\end{equation}
to obtain
\begin{equation}
\frac{\rho }{G_{e}}\,\frac{s^{2}}{1\,+t_{\varepsilon }s}\,%
\tilde{r}\left( x,s\right) =\frac{\partial ^{2}}{\partial x^{2}}\tilde{r}%
\left( x,s\right) .  \label{Eq_mu_Kelvin-Voigt}
\end{equation}%
Substituting a solution of the form
\begin{equation*}
\tilde{r}\left( x,s\right) =a\left( s\right) \exp \left( b\left( s\right)
x\right) ,
\end{equation*}%
we arrive at%
\begin{equation*}
\tilde{r}\left( x,s\right) =a\left( s\right) \exp \left( \pm \sqrt{\frac{%
\rho }{G_{e}}}\,\frac{s}{\sqrt{1\,+t_{\varepsilon }s}}x\right) .
\end{equation*}%
According to the boundary conditions in the Laplace domain (see \cite{MainardiBook}, Sect. 4.2.1), i.e.,
\begin{eqnarray}
\tilde{r}\left( 0,s\right) &=&\tilde{r}_{0}\left( s\right) ,
\label{Cond_boundary_Laplace} \\
\lim_{x\rightarrow \infty }\tilde{r}\left( x,s\right) &=&0,
\label{Cond_boundary_Laplace_limit}
\end{eqnarray}
we have%
\begin{equation*}
\tilde{r}\left( x,s\right) =\tilde{r}_{0}\left( s\right) \exp \left( -\sqrt{%
\frac{\rho }{G_{e}}}\,\frac{s}{\sqrt{1\,+t_{\varepsilon }s}}x\right) .
\end{equation*}%
Since%
\begin{equation*}
\left[ \sqrt{\frac{\rho }{G_{e}}}\right] =\frac{1}{L\,T^{-1}},
\end{equation*}%
has dimensions of the inverse of a velocity, define the velocity
\begin{equation*}
c^{\prime }=\sqrt{\frac{G_{e}}{\rho }},
\end{equation*}%
and recast the solution in the Laplace domain as%
\begin{equation*}
\tilde{r}\left( x,s\right) =\tilde{r}_{0}\left( s\right) \exp \left( -\,%
\frac{s\,}{\sqrt{1\,+t_{\varepsilon }s}}\frac{x}{c^{\prime }}\right) .
\end{equation*}%
Recall that according to the Laplace transform definition \cite{Schiff} (Eqn. 1.1), the dimension of the variable $s$ is $\left[ s\right] =T^{-1}$; thus, $%
t_{\varepsilon }s$ is dimensionless, i.e., $\left[ t_{\varepsilon }s\right]
=1 $.
Apply the convolution theorem of the Laplace transform \cite{Schiff} (Theorem 2.39), to obtain%
\begin{eqnarray}
r\left( x,t\right) &=&\mathcal{L}^{-1}\left[ \tilde{r}_{0}\left( s\right)
\exp \left( -\,\frac{s\,}{\sqrt{1\,+t_{\varepsilon }s}}\frac{x}{c^{\prime }}%
\right);t \right]  \label{u(x,t)_a} \\
&=&\int_{0}^{t}r_{0}\left( t-t^{\prime }\right) \,\mathcal{L}^{-1}\left[
\exp \left( -\,\frac{s\,}{\sqrt{1\,+t_{\varepsilon }s}}\frac{x}{c^{\prime }}%
\right) ;t^{\prime }\right] \,dt^{\prime }.  \notag
\end{eqnarray}%
Now, consider the dimensionless variables:
\begin{eqnarray}
\xi &=&\frac{x}{c^{\prime }\,t_{\varepsilon }},  \label{Xi_def} \\
\tau &=&\frac{t}{t_{\varepsilon }},  \label{Tau_def}
\end{eqnarray}%
and the Laplace transform property \cite{Prudnikov5} (1.1.1(4)): 
\begin{equation}
\mathcal{L}\left[\frac{1}{a}e^{-bt/a}\,f\left( \frac{t}{a}\right);s \right] = \tilde{f}\left( a\,s+b\right)  ,\quad a>0,  \label{L[f]_property}
\end{equation}%
to calculate%
\begin{eqnarray}
\mathcal{L}^{-1}\left[ \exp \left( -\,\frac{s\,}{\sqrt{1\,+t_{\varepsilon }s}%
}\frac{x}{c^{\prime }}\right) ;t\right] &=&\mathcal{L}^{-1}\left[ \exp
\left( -\,\frac{\,\xi \,t_{\varepsilon }s\,}{\sqrt{1\,+t_{\varepsilon }s}}%
\right) ;t\right]  \notag \\
&=&\frac{e^{-\tau }}{t_{\varepsilon }}\mathcal{L}^{-1}\left[ \exp \left(
-\,\xi \sqrt{s}+\frac{\xi }{\sqrt{s}}\right) ;\tau \right]  \label{L-1_trans}
\end{eqnarray}%
In order to calculate the inverse Laplace transform given in (\ref{L-1_trans}%
), define%
\begin{equation}
f\left( \xi ,\tau \right) :=\mathcal{L}^{-1}\left[ \frac{\exp \left( -\,\xi
\sqrt{s}+\xi /\sqrt{s}\,\right) }{s};\tau \right] ,  \label{f(xi,tau)_def}
\end{equation}%
According to the derivative theorem of the Laplace transform \cite{Schiff} (Theorem 2.7), and
knowing that $\mathcal{L}^{-1}\left[ 1\right] =\delta \left( \tau \right) $
\cite{Schiff} (Eqn. 2.37), we have
\begin{equation}
\mathcal{L}^{-1}\left[ \exp \left( -\,\xi \sqrt{s}+\frac{\xi }{\sqrt{s}}%
\right) ;\tau \right] =\frac{\partial }{\partial \tau }f\left( \xi ,\tau
\right) +f\left( \xi ,0\right) \,\delta \left( \tau \right) .
\label{Parcial_f(0)}
\end{equation}%
Now, expand $\exp \left( \xi /\sqrt{s}\,\right) $ in (\ref{f(xi,tau)_def})
to obtain
\begin{equation*}
f\left( \xi ,\tau \right) =\sum_{n=0}^{\infty }\frac{\xi ^{n}}{n!}\mathcal{L}%
^{-1}\left( \frac{\exp \left( -\,\xi \sqrt{s}\right) }{s^{1+n/2}};\tau
\right) .
\end{equation*}%
According to \cite{Prudnikov5} (Eqn. 2.2.1(11)), we have%
\begin{equation*}
\mathcal{L}^{-1}\left( \frac{\exp \left( -\,\xi \sqrt{s}\right) }{s^{n/2+1}}%
;\tau \right) =\frac{2^{n+1}}{\sqrt{\pi }}\exp \left( -\frac{\xi ^{2}}{4\tau
}\right) t^{n/2}\,H_{-1-n}\left( \frac{\xi }{2\sqrt{\tau }}\right) ,
\end{equation*}%
where $H_{\nu }\left( z\right) $ denotes the Hermite function \cite{Lebedev} (Eqn.
10.2.8). Therefore,
\begin{equation}
f\left( \xi ,\tau \right) =\frac{2}{\sqrt{\pi }}\exp \left( -\frac{\xi ^{2}}{%
4\tau }\right) \sum_{n=0}^{\infty }\frac{\left( 2\xi \sqrt{\tau }\right) ^{n}%
}{n!}H_{-1-n}\left( \frac{\xi }{2\sqrt{\tau }}\right) .  \label{f_Hermite}
\end{equation}%
Now, perform the change of variables
\begin{equation*}
z\left( \xi ,\tau \right) =\frac{\xi }{2\sqrt{\tau }},
\end{equation*}%
to obtain%
\begin{equation*}
f\left( \xi ,z\right) =\frac{2}{\sqrt{\pi }}e^{-z^{2}}\sum_{n=0}^{\infty }%
\frac{\left( \xi ^{2}/z\right) ^{n}}{n!}H_{-1-n}\left( z\right) .
\end{equation*}%
According to the integral representation \cite{Lebedev} (Eqn. 10.5.2):%
\begin{equation*}
H_{-\nu }\left( z\right) =\frac{1}{\Gamma \left( \nu \right) }%
\int_{0}^{\infty }\exp \left( -u^{2}-2u\,z\right) \,u^{-\nu -1}du,\quad \Re
\left( \nu \right) >0,
\end{equation*}%
we find that%
\begin{eqnarray*}
f\left( \xi ,z\right) &=&\frac{2}{\sqrt{\pi }}e^{-z^{2}}\sum_{n=0}^{\infty }%
\frac{\left( \xi ^{2}/z\right) ^{n}}{\left( n!\right) ^{2}}\int_{0}^{\infty
}\exp \left( -u^{2}-2u\,z\right) u^{n}\,du \\
&=&\frac{2}{\sqrt{\pi }}\sum_{n=0}^{\infty }\frac{\left( \xi ^{2}/z\right)
^{n}}{\left( n!\right) ^{2}}\int_{0}^{\infty }\exp \left( -\left( u+z\right)
^{2}\right) u^{n}\,du \\
&=&\frac{2}{\sqrt{\pi }}\int_{0}^{\infty }\exp \left( -\left( u+z\right)
^{2}\right) \left[ \sum_{n=0}^{\infty }\left( \frac{\xi ^{2}u}{z}\right) ^{n}%
\frac{1}{\left( n!\right) ^{2}}\right] \,du .
\end{eqnarray*}%
By using the definition of the normalized hypergeometric function $_{0}\overline{F}_{1}$
\cite{Lebedev} (Sect. 9.14), i.e.,
\begin{equation}
_{0}\overline{F}_{1}\left(
\begin{array}{c}
- \\
\gamma%
\end{array}%
;x\right) =\frac{1}{\Gamma \left( \gamma \right) }\,_{0}F_{1}\left(
\begin{array}{c}
- \\
\gamma%
\end{array}%
;x\right) =\frac{1}{\Gamma \left( \gamma \right) }\sum_{k=0}^{\infty }\frac{%
x^{k}}{\left( \gamma \right) _{k}k!},  \label{0F1_regularized_def}
\end{equation}%
we recast the above result as%
\begin{equation}
f\left( \xi ,\tau \right) =\frac{2}{\sqrt{\pi }}\int_{0}^{\infty }\exp
\left( -\left[ u+\frac{\xi }{2\sqrt{\tau }}\right] ^{2}\right) \,_{0}\overline{F}_{1}\left(
\begin{array}{c}
- \\
1%
\end{array}%
;2\xi \sqrt{\tau }u\right) \,du.  \label{f(xi,tau)_resultado}
\end{equation}
Moreover, simplify the above expression performing the change of variables $v=2%
\sqrt{\tau }u-\xi $ to arrive at%
\begin{equation}
f\left( \xi ,\tau \right) =\frac{1}{\sqrt{\pi \tau }}\int_{\xi }^{\infty
}\exp \left( -\frac{v^{2}}{4\tau }\right) \,_{0}\overline{F}_{1}\left(
\begin{array}{c}
- \\
1%
\end{array}%
;\xi \left( v-\xi \right) \right) \,dv.  \label{f(xi,tau)_resultado2}
\end{equation}
Note that
\begin{equation}
\lim_{\tau \rightarrow 0}f\left( \xi ,\tau \right) =0,  \label{f(xi,0)=0}
\end{equation}%
Thus, according to (\ref{Parcial_f(0)}), we have
\begin{equation}
g\left( \xi ,\tau \right) :=\frac{\partial }{\partial \tau }f\left( \xi
,\tau \right) =\mathcal{L}^{-1}\left[ \exp \left( -\,\xi \sqrt{s}+\frac{\xi
}{\sqrt{s}}\right) ;\tau \right] .  \label{g(xi,tau)_def}
\end{equation}%
However, from (\ref{f(xi,tau)_resultado2}), we arrive at%
\begin{equation}
g\left( \xi ,\tau \right) =\frac{1}{4\sqrt{\pi }\,\tau ^{5/2}}\int_{\xi
}^{\infty }\left( v^{2}-2\tau \right) \exp \left( -\frac{v^{2}}{4\tau }%
\right) \,_{0}\overline{F}_{1}\left(
\begin{array}{c}
- \\
1%
\end{array}%
;\xi \left( v-\xi \right) \right) dv.  \label{g(xi,tau)_resultado}
\end{equation}%
Consequently,
\begin{equation*}
\mathcal{L}^{-1}\left[ \exp \left( -\,\frac{\,\xi \,t_{\varepsilon }s\,}{%
\sqrt{1\,+t_{\varepsilon }s}}\right) ;t\right] =\frac{\exp \left(
-t/t_{\varepsilon }\right) }{t_{\varepsilon }}g\left( \xi ,t/t_{\varepsilon
}\right) ,
\end{equation*}%
and
\begin{equation}
r\left( \xi ,t\right) =\frac{1}{t_{\varepsilon }}\int_{0}^{t}r_{0}\left(
t-t^{\prime }\right) \,\exp \left( -\frac{t^{\prime }}{t_{\varepsilon }}%
\right) \,g\left( \xi ,\frac{t^{\prime }}{t_{\varepsilon }}\right)
\,dt^{\prime }.  \label{u(x,t)_general}
\end{equation}%
Now, perform the change of variables $\tau ^{\prime }=t^{\prime
}/t_{\varepsilon }$ and take into account the dimensionless variable $\tau
=t/t_{\varepsilon }$ given in (\ref{Tau_def}) to obtain%
\begin{equation}
r\left( \xi ,\tau \right) =\int_{0}^{\tau }r_{0}\left( t_{\varepsilon
}\left( \tau -\tau ^{\prime }\right) \right) \,e^{-\tau ^{\prime }}\,g\left(
\xi ,\tau ^{\prime }\right) \,d\tau ^{\prime }.
\label{u(x,t)_general_dimensionless}
\end{equation}

\subsection{Step Pulse Solution}

For the case $r_{0}\left( t\right) =\theta \left( t\right) $, where $\theta \left( t\right)$ denotes the Heaviside function, we state that
\begin{equation*}
\tilde{r}_{0}\left( s\right) =\frac{1}{s},
\end{equation*}%
hence, (\ref{u(x,t)_a})\ results in
\begin{equation*}
r\left( x,t\right) =\mathcal{L}^{-1}\left[ \frac{1}{s}\exp \left( -\,\frac{%
s\,}{\sqrt{1\,+t_{\varepsilon }s}}\frac{x}{c^{\prime }}\right) ;t\right] .
\end{equation*}%
According to the dimensionless variable $\xi $ given in (\ref{Xi_def}),
\begin{eqnarray*}
r\left( \xi ,t\right) &=&\mathcal{L}^{-1}\left[ \frac{1}{s}\exp \left( -\,%
\frac{\xi \,t_{\varepsilon }\,s\,}{\sqrt{1\,+t_{\varepsilon }s}}\right) ;t%
\right] \\
&=&t_{\varepsilon }\,\mathcal{L}^{-1}\left[ \frac{1}{t_{\varepsilon }\,s}%
\exp \left( -\,\frac{\xi \,t_{\varepsilon }\,s\,}{\sqrt{1\,+t_{\varepsilon }s%
}}\right) ;t\right] ,
\end{eqnarray*}%
and applying (\ref{Tau_def})\ and (\ref{L[f]_property}), we get%
\begin{equation}
r\left( \xi ,\tau \right) =\mathcal{L}^{-1}\left[ \frac{1}{s}\exp \left( -\,%
\frac{\xi \,s\,}{\sqrt{1\,+s}}\right) ;\tau \right] ,
\label{u(xi,tau)_Laplace}
\end{equation}%
so, as $s\rightarrow 0$, i.e., $\tau \rightarrow \infty $ (or equivalently, $%
\xi \rightarrow 0$) we arrive at the following asymptotic representation:\
\begin{equation}
r\left( \xi ,\tau \right) \approx \mathcal{L}^{-1}\left[ \frac{1}{s};\tau %
\right] =1,\quad \tau \rightarrow \infty ,\quad \left( \xi \rightarrow
0\right) .  \label{u(xi,tau)_tau->inf}
\end{equation}
However, if we take $r_{0}\left( t\right) =\theta \left( t\right) $ in (\ref{u(x,t)_general_dimensionless}), we obtain
\begin{equation*}
r\left( \xi ,\tau \right) =\int_{0}^{\tau }e^{-\tau ^{\prime }}g\left( \xi
,\tau ^{\prime }\right) \,d\tau ^{\prime }.
\end{equation*}%
Integrating by parts, taking into account (\ref{f(xi,0)=0}), we arrive at%
\begin{equation}
r\left( \xi ,\tau \right) =e^{-\tau }f\left( \xi ,\tau \right)
+\int_{0}^{\tau }e^{-\tau ^{\prime }}f\left( \xi ,\tau ^{\prime }\right)
\,d\tau ^{\prime },  \label{u(xi,tau)_theta}
\end{equation}%
where $f\left( \xi ,\tau \right) $ is given in integral form in (\ref{f(xi,tau)_resultado2}). According to the integral \cite{Prudnikov1} (Eqn. 1.3.3(20))
\begin{equation*}
\int_{0}^{\tau }\frac{\exp \left( -\tau ^{\prime }-\frac{v^{2}}{4\tau
^{\prime }}\right) }{\sqrt{\pi \,\tau ^{\prime }}}\,d\tau ^{\prime }=\frac{1%
}{2}\left[ e^{-v}\,\mathrm{erfc}\left( \frac{v-2\tau }{2\sqrt{\tau }}\right)
-e^{v}\,\mathrm{erfc}\left( \frac{v+2\tau }{2\sqrt{\tau }}\right) \right] ,
\end{equation*}%
the solution given in (\ref{u(xi,tau)_theta})\ is reduced to
\begin{eqnarray}
r\left( \xi ,\tau \right)  &=&\int_{\xi }^{\infty }\,_{0}\overline{F}_{1}\left(
\begin{array}{c}
- \\
1%
\end{array}%
;\xi \left( v-\xi \right) \right)   \label{u(xi,tau)_theta_2} \\
&&\left[ \frac{\exp \left( -\frac{v^{2}}{4\tau }-\tau \right) }{\sqrt{\pi
\,\tau }}+\frac{e^{-v}}{2}\mathrm{erfc}\left( \frac{v-2\tau }{2\sqrt{\tau }}%
\right) -\frac{e^{v}}{2}\,\mathrm{erfc}\left( \frac{v+2\tau }{2\sqrt{\tau }}%
\right) \right] dv.  \notag
\end{eqnarray}

\subsubsection{Asymptotic Solution for $\protect\tau \rightarrow 0$ or $\protect\xi \rightarrow \infty $}

Recall that $f\left( \xi ,\tau \right) $ is given by (\ref{f(xi,tau)_resultado}), thus
\begin{eqnarray}
f\left( \xi ,\tau \right) &\approx &\frac{2}{\sqrt{\pi }}\int_{0}^{\infty
}\exp \left( -\left[ u+\frac{\xi }{2\sqrt{\tau }}\right] ^{2}\right) \,\,du
\notag \\
&=&\mathrm{erfc}\left( \frac{\xi }{2\sqrt{\tau }}\right) ,\quad \tau
\rightarrow 0 \quad \text{or}\quad \xi \rightarrow \infty .  \label{f(xi,tau)_ta->0}
\end{eqnarray}%
Therefore, from (\ref{u(xi,tau)_theta})\ and (\ref{f(xi,tau)_ta->0}), we
arrive at the following asymptotic representation for $\tau \rightarrow 0$ (for fixed $\xi$) or $\xi \rightarrow \infty$ (for fixed $\tau$):
\begin{equation}
r\left( \xi ,\tau \right) \approx \mathrm{erfc}\left( \frac{\xi }{2\sqrt{%
\tau }}\right) ,\quad \tau \rightarrow 0 \quad \text{or}\quad \xi \rightarrow \infty .
\label{u(xi,tau)_tau->0}
\end{equation}

\subsubsection{Asymptotic Solution for $\protect\tau \rightarrow \infty $%
}

Rewrite (\ref{u(xi,tau)_theta})\ as
\begin{equation}
r\left( \xi ,\tau \right) =e^{-\tau }f\left( \xi ,\tau \right) +\underset{I}{%
\underbrace{\int_{0}^{\infty }e^{-t}f\left( \xi ,t\right) \,dt}}-\underset{%
I\left( \xi ,\tau \right) }{\underbrace{\int_{\tau }^{\infty }e^{-t}f\left(
\xi ,t\right) \,dt}}.  \label{u_split}
\end{equation}
First, we calculate
\begin{equation}
I:=\int_{0}^{\infty }e^{-t}f\left( \xi ,t\right) \,dt.
\label{I_infinity_def}
\end{equation}%
For this purpose, expand $f\left( \xi ,\tau \right) $ according to (\ref{0F1_regularized_def})\ and (\ref{f(xi,tau)_resultado}) as
\begin{equation*}
f\left( \xi ,\tau \right) =\frac{2}{\sqrt{\pi }}\sum_{k=0}^{\infty }\frac{%
\left( 2\xi \sqrt{\tau }\right) ^{k}}{\left( k!\right) ^{2}}\int_{0}^{\infty
}\exp \left( -\left[ u+\frac{\xi }{2\sqrt{\tau }}\right] ^{2}\right)
\,u^{k}\,du,
\end{equation*}%
where MATHEMATICA's computer algebra yields the following result (recall the definition of the $_{p}F_{q}(z)$ generalized hypergeometric function \cite{NIST} (Sect. 16.2))
\begin{eqnarray*}
&&\int_{0}^{\infty }\exp \left( -\left[ u+\frac{\xi }{2\sqrt{\tau }}\right]
^{2}\right) \,u^{k}\,du \\
&=&\frac{1}{2}\left\{ -\frac{\xi }{\sqrt{\tau }}\Gamma \left( 1+\frac{k}{2}%
\right) \,_{1}F_{1}\left(
\begin{array}{c}
\frac{1-k}{2} \\
\frac{3}{2}%
\end{array}%
;-\frac{\xi ^{2}}{4\tau }\right) +\Gamma \left( \frac{1+k}{2}\right)
\,_{1}F_{1}\left(
\begin{array}{c}
-\frac{k}{2} \\
\frac{1}{2}%
\end{array}%
;-\frac{\xi ^{2}}{4\tau }\right) \right\} ,
\end{eqnarray*}%
thus%
\begin{eqnarray*}
f\left( \xi ,\tau \right) &=&-\frac{\xi }{\sqrt{\pi }}\sum_{k=0}^{\infty }%
\frac{\left( 2\xi \right) ^{k}}{\left( k!\right) ^{2}}\,\Gamma \left( 1+%
\frac{k}{2}\right) \tau ^{\left( k-1\right) /2}\,_{1}F_{1}\left(
\begin{array}{c}
\frac{1-k}{2} \\
\frac{3}{2}%
\end{array}%
;-\frac{\xi ^{2}}{4\tau }\right) \\
&&+\frac{1}{\sqrt{\pi }}\sum_{k=0}^{\infty }\frac{\left( 2\xi \right) ^{k}}{%
\left( k!\right) ^{2}}\,\Gamma \left( \frac{1+k}{2}\right) \,\tau
^{k/2}\,_{1}F_{1}\left(
\begin{array}{c}
-\frac{k}{2} \\
\frac{1}{2}%
\end{array}%
;-\frac{\xi ^{2}}{4\tau }\right) ,
\end{eqnarray*}%
and%
\begin{eqnarray}
I&=&-\frac{\xi }{\sqrt{\pi }}\sum_{k=0}^{\infty }\frac{\left( 2\xi \right)
^{k}}{\left( k!\right) ^{2}}\,\Gamma \left( 1+\frac{k}{2}\right)
\int_{0}^{\infty }e^{-t}\,t^{\left( k-1\right) /2}\,_{1}F_{1}\left(
\begin{array}{c}
\frac{1-k}{2} \\
\frac{3}{2}%
\end{array}%
;-\frac{\xi ^{2}}{4t}\right) dt  \notag \\
&&+\frac{1}{\sqrt{\pi }}\sum_{k=0}^{\infty }\frac{\left( 2\xi \right) ^{k}}{%
\left( k!\right) ^{2}}\,\Gamma \left( \frac{1+k}{2}\right) \int_{0}^{\infty
}e^{-t}\,t^{k/2}\,_{1}F_{1}\left(
\begin{array}{c}
-\frac{k}{2} \\
\frac{1}{2}%
\end{array}%
;-\frac{\xi ^{2}}{4t}\right) dt.  \notag
\end{eqnarray}%
The integrals given above can also be
calculated using MATHEMATICA's computer algebra\ as follows:%
\begin{eqnarray*}
&&\int_{0}^{\infty }e^{-t}\,t^{\left( k-1\right) /2}\,_{1}F_{1}\left(
\begin{array}{c}
\frac{1-k}{2} \\
\frac{3}{2}%
\end{array}%
;-\frac{\xi ^{2}}{4t}\right) dt \\
&=&\frac{2^{-2\left( 1+k\right) }\sqrt{\pi }\,k!}{\Gamma \left( 1+\frac{k}{2}%
\right) \xi }\left[ -\sqrt{\pi }\,\xi ^{2+k}\,_{1}\overline{F}_{2}\left(
\begin{array}{c}
1 \\
\frac{3+k}{2},\frac{4+k}{2}%
\end{array}%
;\frac{\xi ^{2}}{4}\right) +2^{2+k}\sinh \xi \right] , \\
&&\int_{0}^{\infty }e^{-t}\,t^{k/2}\,_{1}F_{1}\left(
\begin{array}{c}
-\frac{k}{2} \\
\frac{1}{2}%
\end{array}%
;-\frac{\xi ^{2}}{4t}\right) dt \\
&=&\frac{2^{-2\left( 1+k\right) }\sqrt{\pi }\,k!}{\Gamma \left( \frac{1+k}{2}%
\right) }\left[ 2^{2+k}\cosh \xi -\sqrt{\pi }\,\xi ^{2+k}\,_{1}\overline{F}%
_{2}\left(
\begin{array}{c}
1 \\
\frac{3+k}{2},\frac{4+k}{2}%
\end{array}%
;\frac{\xi ^{2}}{4}\right) \right] ,
\end{eqnarray*}%
thus 
\begin{equation}
I=\left( -\sinh \xi +\cosh \xi \right) \sum_{k=0}^{\infty }\frac{\xi ^{k}}{k!%
}=1,  \label{I=1}
\end{equation}%
which is consistent with the asymptotic behavior given in (\ref{u(xi,tau)_tau->inf}).

Second, we calculate
\begin{equation}
I\left( \xi ,\tau \right) :=\int_{\tau }^{\infty }e^{-t}f\left( \xi
,t\right) \,dt.  \label{I(xi,tau)_def}
\end{equation}%
For this purpose, apply the following expansion of the Hermite functions
\cite{Lebedev} (10.4.3):%
\begin{equation*}
H_{\nu }\left( z\right) =\frac{1}{2\,\Gamma \left( -\nu \right) }%
\sum_{m=0}^{\infty }\frac{\left( -1\right) ^{m}\Gamma \left( \frac{m-\nu }{2}%
\right) }{m!}\left( 2z\right) ^{m},\quad \left\vert z\right\vert <\infty ,
\end{equation*}%
to obtain%
\begin{equation*}
H_{-1-n}\left( z\right) =\frac{1}{2n!}\left[ \Gamma \left( \frac{1+n}{2}%
\right) -\Gamma \left( 1+\frac{n}{2}\right) 2z+\cdots \right] ,
\end{equation*}%
thus, taking into account the expansion given in (\ref{f_Hermite}), we have
the asymptotic expansion as $t\rightarrow \infty $,%
\begin{eqnarray}
f\left( \xi ,t\right)
\approx \frac{2}{\sqrt{\pi }}\left[ \sum_{n=0}^{\infty }\frac{\left( 2\xi
\sqrt{t}\right) ^{n}}{2\left( n!\right) ^{2}}\Gamma \left( \frac{1+n}{2}%
\right) -\frac{\xi }{t}\sum_{n=0}^{\infty }\frac{\left( 2\xi \sqrt{t}\right)
^{n}}{2\left( n!\right) ^{2}}\Gamma \left( 1+\frac{n}{2}\right) \right] .
\label{f_Hermite_expansion}
\end{eqnarray}%
Now, apply the duplication formula \cite{Lebedev} (Eqn. 1.2.3)%
\begin{equation*}
\sqrt{\pi }\,\Gamma \left( 2z\right) =2^{2z-1}\,\Gamma \left( z\right)
\,\Gamma \left( z+\frac{1}{2}\right) ,
\end{equation*}%
to recast (\ref{f_Hermite_expansion})\ as%
\begin{equation}
f\left( \xi ,t\right) \approx \sum_{n=0}^{\infty }\frac{\left( \xi \sqrt{t}%
\right) ^{n}}{n!\,\Gamma \left( 1+\frac{n}{2}\right) }-\frac{\xi }{\sqrt{t}}%
\sum_{n=0}^{\infty }\frac{\left( \xi \sqrt{t}\right) ^{n}}{n!\,\Gamma \left(
\frac{1+n}{2}\right) },\quad t\rightarrow \infty ,  \label{f_expansion_1}
\end{equation}%
where we have also,
\begin{equation}
f\left( \xi ,t\right) \approx \sum_{n=0}^{\infty }\frac{\left( \xi \sqrt{t}%
\right) ^{n}}{n!\,\Gamma \left( 1+\frac{n}{2}\right) },\quad t\rightarrow
\infty .  \label{f_expansion_2}
\end{equation}%
Now, insert (\ref{f_expansion_2})\ into (\ref{I(xi,tau)_def})\ to obtain%
\begin{equation*}
I\left( \xi ,\tau \right) \approx \sum_{n=0}^{\infty }\frac{\xi ^{n}}{%
n!\,\Gamma \left( 1+\frac{n}{2}\right) }\Gamma \left( 1+\frac{n}{2},\tau
\right) ,\quad \tau \rightarrow \infty ,
\end{equation*}%
where we have used the upper incomplete gamma function \cite{Atlas} (Eqn. 45:3:2):%
\begin{equation*}
\Gamma \left( \nu ,z\right) :=\int_{z}^{\infty }e^{-t}\,t^{\nu -1}\,dt.
\end{equation*}%
Using the asymptotic expansion \cite{NIST} (Eqns. 8.11.1-3):%
\begin{equation*}
\Gamma \left( \nu ,z\right) \approx z^{\nu -1}e^{-z}\left(
\sum_{k=0}^{\infty }\frac{\left( -1\right) ^{k}\left( 1-\nu \right) _{k}}{%
z^{k}}+O\left( z^{-n}\right) \right) ,\quad z\rightarrow \infty ,
\end{equation*}%
we arrive at%
\begin{equation}
I\left( \xi ,\tau \right) \approx e^{-\tau }\left[ \sum_{n=0}^{\infty }\frac{%
\left( \xi \sqrt{\tau }\right) ^{n}}{n!\,\Gamma \left( 1+\frac{n}{2}\right) }%
-\frac{1}{2\tau }\sum_{n=0}^{\infty }\frac{\left( \xi \sqrt{\tau }\right)
^{n}\,n}{n!\,\Gamma \left( 1+\frac{n}{2}\right) }\right] ,\quad \tau
\rightarrow \infty .  \label{I(xi,tau) sim1}
\end{equation}%
Recast the second sum in (\ref{I(xi,tau) sim1}), to obtain,%
\begin{equation}
I\left( \xi ,\tau \right) \approx e^{-\tau }\left[ \sum_{n=0}^{\infty }\frac{%
\left( \xi \sqrt{\tau }\right) ^{n}}{n!\,\Gamma \left( 1+\frac{n}{2}\right) }%
-\frac{\xi }{2\sqrt{\tau }}\sum_{n=0}^{\infty }\frac{\left( \xi \sqrt{\tau }%
\right) ^{n}\,}{n!\,\Gamma \left( \frac{3+n}{2}\right) }\right] ,\quad \tau
\rightarrow \infty .  \label{I(xi,tau)_tau->inf}
\end{equation}

Finally, insert\ the results given in (\ref{I=1}), (\ref{f_expansion_1})\
and (\ref{I(xi,tau)_tau->inf}), in (\ref{u_split}), and simplify the result \mbox{to obtain}
\begin{equation*}
r\left( \xi ,\tau \right) \approx 1-\frac{\xi \,e^{-\tau }}{\sqrt{\tau }}%
\sum_{n=0}^{\infty }\frac{\left( \xi \sqrt{\tau }\right) ^{n}}{n!\,}\left[
\frac{1}{\Gamma \left( 1+\frac{n}{2}\right) }-\frac{1}{2\,\Gamma \left(
\frac{3+n}{2}\right) }\right] ,\quad \tau \rightarrow \infty .
\end{equation*}%
Apply the factorial property of the gamma function \cite{Lebedev} (Eqn. 1.2.1)%
, i.e., $\Gamma \left( z+1\right) =z\,\Gamma \left( z\right) $, to
reformulate the above sum as%
\begin{equation*}
r\left( \xi ,\tau \right) \approx 1-\xi ^{2}\,e^{-\tau }\sum_{n=0}^{\infty }%
\frac{\left( \xi \sqrt{\tau }\right) ^{n}}{n!\,\left( n+2\right) \,\Gamma
\left( 1+\frac{n}{2}\right) },\quad \tau \rightarrow \infty .
\end{equation*}%
Split the above sum in the even and odd terms and recast the result as two
hypergeometric sums to arrive at the following asymptotic expression,
\begin{eqnarray}
&&r\left( \xi ,\tau \right)
\approx 1-\xi ^{2}\,e^{-\tau }\left[ \frac{1}{2}\,_{0}F_{2}\left(
\begin{array}{c}
- \\
\frac{1}{2},2%
\end{array}%
;\frac{\xi ^{2}\tau }{4}\right) +\frac{2\xi \sqrt{\tau }}{3\sqrt{\pi }}%
\,_{0}F_{2}\left(
\begin{array}{c}
- \\
\frac{3}{2},\frac{5}{2}%
\end{array}%
;\frac{\xi ^{2}\tau }{4}\right) \right] ,  \notag \\
&&\tau \rightarrow \infty .  \label{u(xi,tau)_tau_inf2b}
\end{eqnarray}

\begin{figure}[htbp]
\includegraphics[width=0.8 \textwidth]{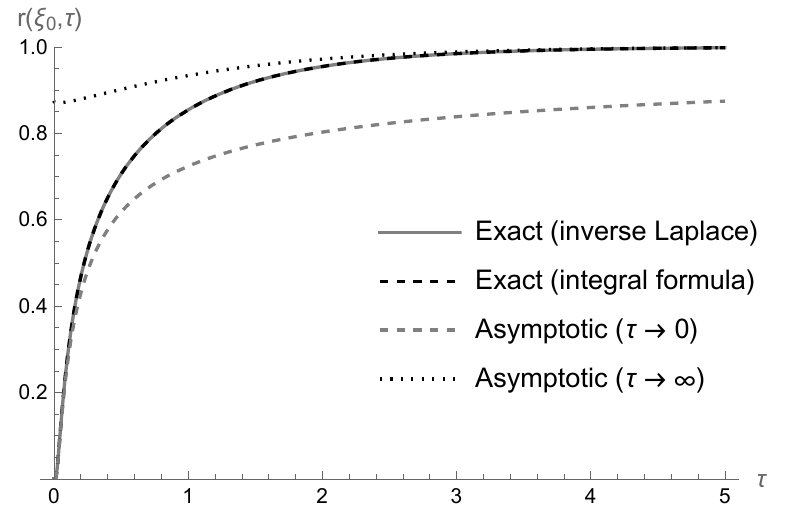}
\caption{Graph of $r\left( \protect\xi _{0},\protect\tau \right) $ with $%
\protect\xi _{0}=0.5$ for the step pulse solutions and the asymptotic
approximations in the Kelvin--Voigt model. }
\label{Figure: Kelvin-Voigt}
\end{figure}

Figure \ref{Figure: Kelvin-Voigt} shows that the solution $r\left( \xi ,\tau
\right) $ given in terms of the inverse Laplace transform~(\ref{u(xi,tau)_Laplace})\ is numerically equivalent to the integral solution (\ref{u(xi,tau)_theta_2}).
Indeed, the maximum discrepancy between both expressions in
Figure \ref{Figure: Kelvin-Voigt} is $8.30117\times 10^{-10}$, thus
within the numerical error of (\ref{u(xi,tau)_Laplace}) and (\ref{u(xi,tau)_theta_2}) computed by MATHEMATICA
(see the MATHEMATICA notebook available at 
 \url{https://shorturl.at/d0bgQ} accessed on 11 January 2026).
Also, the asymptotic representations for $\tau
\rightarrow \infty $, i.e., (\ref{u(xi,tau)_tau_inf2b}),\ and $\tau
\rightarrow 0$, i.e., (\ref{u(xi,tau)_tau->0}) are numerically checked.

\subsubsection{Asymptotic Solution for $\protect\xi \rightarrow 0$}

In order to obtain the asymptotic solution as $\xi \rightarrow 0$, expand
the integrand of $f\left( \xi ,\tau \right) $ in powers of $\xi $ as
follows:
\begin{eqnarray*}
\exp \left( -\left[ u+\frac{\xi }{2\sqrt{\tau }}\right] ^{2}\right)
&=&e^{-u^{2}}\exp \left( -\frac{\xi }{\sqrt{\tau }}\left( u+\frac{\xi }{4%
\sqrt{\tau }}\right) \right) \\
&=&e^{-u^{2}}\sum_{k=0}^{\infty }\frac{\left( -1\right) ^{k}}{k!}\left(
\frac{\xi }{\sqrt{\tau }}\right) ^{k}\left( u+\frac{\xi }{4\sqrt{\tau }}%
\right) ^{k} \\
&=&e^{-u^{2}}\left[ 1-\frac{\xi }{\sqrt{\tau }}\left( u+\frac{\xi }{4\sqrt{%
\tau }}\right) +\frac{\xi ^{2}}{2\tau }\left( u+\frac{\xi }{4\sqrt{\tau }}%
\right) ^{2}+\cdots \right] ,
\end{eqnarray*}%
and%
\begin{eqnarray*}
_{0}\overline{F}_{1}\left(
\begin{array}{c}
- \\
1%
\end{array}%
;2\xi \sqrt{\tau }u\right) &=&\sum_{k=0}^{\infty }\frac{\left( 2\xi \sqrt{%
\tau }u\right) ^{k}}{\left( k!\right) ^{2}} \\
&=&1+2\xi \sqrt{\tau }u+\xi ^{2}\,\tau \,u^{2}+\cdots
\end{eqnarray*}%
Therefore, up to the second order, the integrand of (\ref{f(xi,tau)_resultado})\
is given by
\begin{eqnarray}
&&\exp \left( -\left[ u+\frac{\xi }{2\sqrt{\tau }}\right] ^{2}\right) \,_{0}%
\overline{F}_{1}\left(
\begin{array}{c}
- \\
1%
\end{array}%
;2\xi \sqrt{\tau }u\right)  \label{integrand_e->0} \\
&\approx &e^{-u^{2}}\left\{ 1+\xi \frac{u}{\sqrt{\tau }}\left( 2\tau
-1\right) +\xi ^{2}\left[ u^{2}\left( \frac{1}{2\tau }+\tau -2\right) -\frac{%
1}{4\tau }\right] \right\} ,\quad \xi \rightarrow 0.  \notag
\end{eqnarray}%
Insert (\ref{integrand_e->0})\ into (\ref{f(xi,tau)_resultado}) and apply
the integral%
\begin{equation*}
\int_{0}^{\infty }e^{-u^{2}}u^{k}du=\frac{1}{2}\,\Gamma \left( \frac{1+k}{2}%
\right) ,
\end{equation*}%
to arrive at%
\begin{equation}
f\left( \xi ,\tau \right) \approx 1+\xi \,\frac{2\tau -1}{\sqrt{\pi \,\tau }}%
+\frac{\xi ^{2}}{2}\left( \tau -2\right) ,\quad \xi \rightarrow 0,
\label{f(xi,tau),xi->0}
\end{equation}%
which is consistent with (\ref{u(xi,tau)_tau->inf}). Moreover, taking into
account (\ref{f(xi,tau),xi->0}), we have
\begin{eqnarray}
&&\int_{0}^{\tau }e^{-\tau ^{\prime }}f\left( \xi ,\tau ^{\prime }\right)
\,d\tau ^{\prime }  \label{int(exp(-t)*f)_xi->0} \\
&\approx &\int_{0}^{\tau }e^{-\tau ^{\prime }}\left[ 1+\xi \,\frac{2\tau
^{\prime }-1}{\sqrt{\pi \,\tau ^{\prime }}}+\frac{\xi ^{2}}{2}\left( \tau
^{\prime }-2\right) \right] \,d\tau ^{\prime }  \notag \\
&=&1-e^{-\tau }-\frac{2\sqrt{\tau }e^{-\tau }}{\sqrt{\pi }}\,\xi +\frac{%
\left( 1-\tau \right) e^{-\tau }-1}{2}\,\xi ^{2},\quad \xi \rightarrow 0.
\notag
\end{eqnarray}%
Finally, insert (\ref{f(xi,tau),xi->0})\ and (\ref{int(exp(-t)*f)_xi->0}) in
(\ref{u(xi,tau)_theta})\ and simplify the result to arrive at%
\begin{equation}
r\left( \xi ,\tau \right) \approx 1-\frac{e^{-\tau }}{\sqrt{\pi \,\tau }}%
\,\xi -\frac{1+e^{-\tau }}{2}\,\xi ^{2},\quad \xi \rightarrow 0.
\label{u(xi,tau)_xi->0}
\end{equation}

\begin{figure}[htbp]
\includegraphics[width=0.8
\textwidth]{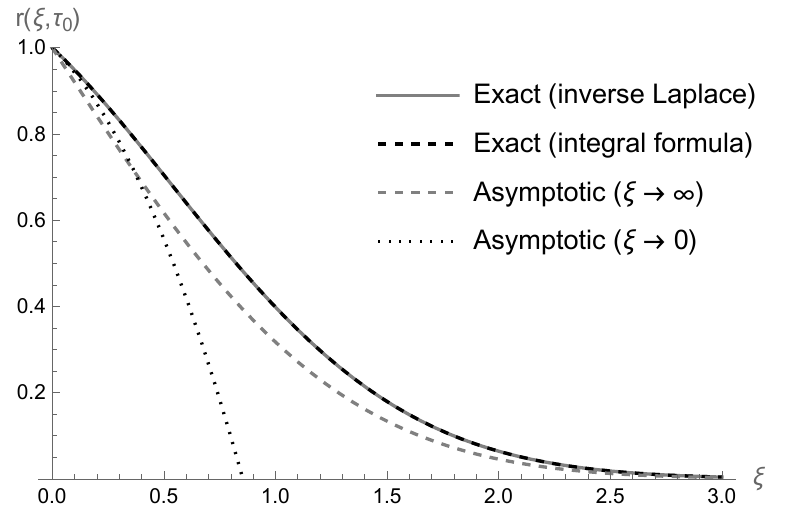}
\caption{Graph of $r\left( \protect\xi ,\protect\tau _{0}\right) $ with $%
\protect\tau _{0}=0.5$ for the step pulse solutions and the asymptotic
approximations in the Kelvin--Voigt model. }
\label{Figure: Kelvin-Voigt_xi}
\end{figure}

Figure \ref{Figure: Kelvin-Voigt_xi} shows that the solution $r\left( \xi
,\tau \right) $ given in terms of the inverse Laplace transform~(\ref{u(xi,tau)_Laplace})\ is numerically equivalent to the integral solution (\ref{u(xi,tau)_theta_2}).
Indeed, the maximum discrepancy between both expressions in
Figure \ref{Figure: Kelvin-Voigt_xi} is $8.35081\times10^{-11}$, thus within the numerical error
of (\ref{u(xi,tau)_Laplace}) and (\ref{u(xi,tau)_theta_2}) computed by MATHEMATICA
(see the MATHEMATICA notebook available at
 \url{https://shorturl.at/d0bgQ} accessed on 11 January 2026).
Also, the asymptotic representations for $\xi
\rightarrow \infty $, i.e., (\ref{u(xi,tau)_tau->0}), and $\xi \rightarrow 0$%
, i.e., (\ref{u(xi,tau)_xi->0}), are numerically checked. Notice as well
that, according to~(\ref{u(xi,tau)_tau->inf}), we have
\begin{equation*}
\lim_{\xi \rightarrow 0}r\left( \xi ,\tau \right) =1=\lim_{\xi \rightarrow 0}%
\mathrm{erfc}\left( \frac{\xi }{2\sqrt{\tau }}\right) ,
\end{equation*}%
as we can see in Figure \ref{Figure: Kelvin-Voigt_xi}.

\subsection{Delta Pulse Solution}

For the case $r_{0}\left( t\right) =\delta \left( t/t_{\varepsilon }\right) $%
, the solution given in (\ref{u(x,t)_general_dimensionless})\ is reduced to%
\begin{eqnarray}
r\left( \xi ,\tau \right)  &=&\int_{0}^{\tau }\delta \left( \tau -\tau
^{\prime }\right) \,e^{-\tau ^{\prime }}\,g\left( \xi ,\tau ^{\prime
}\right) \,d\tau ^{\prime } \notag \\
&=&e^{-\tau }\,g\left( \xi ,\tau \right) , \label{r_delta_resultado}
\end{eqnarray}%
Therefore, according to (\ref{g(xi,tau)_def}), we have%
\begin{eqnarray}
r\left( \xi ,\tau \right)  &=&e^{-\tau }\,\frac{\partial }{\partial \tau }%
f\left( \xi ,\tau \right)   \label{u_delta_derivative_f} \\
&=&e^{-\tau }\,\mathcal{L}^{-1}\left[ \exp \left( -\,\xi \sqrt{s}+\frac{\xi
}{\sqrt{s}}\right) ;\tau \right] ,  \label{u_delta_inv_Laplace}
\end{eqnarray}%
and according to the integral form of $g\left( \xi ,\tau \right) $ given in (%
\ref{g(xi,tau)_resultado}), we have%
\begin{equation}
r\left( \xi ,\tau \right) =\frac{e^{-\tau }}{4\sqrt{\pi }\,\tau ^{5/2}}%
\int_{\xi }^{\infty }\left( v^{2}-2\tau \right) \exp \left( -\frac{v^{2}}{%
4\tau }\right) \,_{0}\overline{F}_{1}\left(
\begin{array}{c}
- \\
1%
\end{array}%
;\xi \left( v-\xi \right) \right) dv.  \label{u_delta_integral}
\end{equation}%
Now, from (\ref{f(xi,tau)_ta->0}) and (\ref{u_delta_derivative_f}), we
arrive at the following asymptotic formula:%
\begin{equation}
r\left( \xi ,\tau \right) \approx \frac{\xi }{2\sqrt{\pi }\,\tau ^{3/2}}\exp
\left( -\frac{\xi ^{2}}{4\tau }-\tau \right) ,\quad \tau \rightarrow 0 \quad \text{or} \quad
\xi \rightarrow \infty .  \label{u_delta_tau->0}
\end{equation}%
Also, recasting (\ref{f_expansion_2}) in terms of generalized hypergeometric functions,
we have%
\begin{equation*}
f\left( \xi ,\tau \right) \approx \,_{0}F_{2}\left(
\begin{array}{c}
- \\
\frac{1}{2},1%
\end{array}%
;\frac{\xi ^{2}\tau }{4}\right) +\frac{2\xi \sqrt{\tau }}{\sqrt{\pi }}%
\,_{0}F_{2}\left(
\begin{array}{c}
- \\
\frac{3}{2},\frac{3}{2}%
\end{array}%
;\frac{\xi ^{2}\tau }{4}\right) ,\quad \tau \rightarrow \infty ,
\end{equation*}%
thus%
\begin{eqnarray}
&&r\left( \xi ,\tau \right) \approx \xi \,e^{-\tau }
\label{u_delta_tau->inf} \\
&&\times\left[ \frac{2\xi ^{2}\sqrt{\tau }}{9\sqrt{\pi }}\,_{0}F_{2}\left(
\begin{array}{c}
- \\
\frac{5}{2},\frac{5}{2}%
\end{array}%
;\frac{\xi ^{2}\tau }{4}\right) +\frac{\xi }{2}\,_{0}F_{2}\left(
\begin{array}{c}
- \\
\frac{3}{2},2%
\end{array}%
;\frac{\xi ^{2}\tau }{4}\right) +\frac{1}{\sqrt{\pi \,\tau }}%
\,_{0}F_{2}\left(
\begin{array}{c}
- \\
\frac{3}{2},\frac{3}{2}%
\end{array}%
;\frac{\xi ^{2}\tau }{4}\right) \right] ,  \notag \\
&&\tau \rightarrow \infty   \notag
\end{eqnarray}%
Finally, from (\ref{f(xi,tau),xi->0}) we have
\begin{equation}
r\left( \xi ,\tau \right) \approx \xi \,e^{-\tau }\left[ \frac{\xi }{2}+%
\frac{1}{\sqrt{\pi \,\tau }}\left( 1+\frac{1}{2\tau }\right) \right] ,\quad
\xi \rightarrow 0. \label{u_delta_xi->0}
\end{equation}%

\begin{figure}[htbp]
\includegraphics[width=0.8 \textwidth]{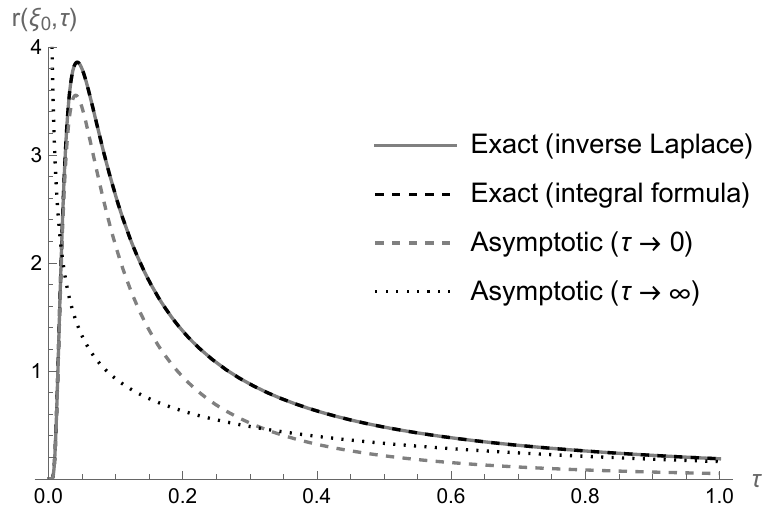}
\caption{Graph of $r\left( \protect\xi _{0},\protect\tau \right) $ with $%
\protect\xi _{0}=0.5$ for the delta pulse solutions and the asymptotic
approximations in the Kelvin--Voigt model. }
\label{Figure: u_delta_tau}
\end{figure}

Figure \ref{Figure: u_delta_tau} shows that the solution $r\left( \xi ,\tau
\right) $ given in terms of the inverse Laplace transform~(\ref{u_delta_inv_Laplace})\ is numerically equivalent to the integral solution (%
\ref{u_delta_integral}).
Indeed, the maximum discrepancy between both expressions in
Figure \ref{Figure: u_delta_tau} is $5.10461\times10^{-9}$, thus within the numerical error
of (\ref{u_delta_inv_Laplace}) and (\ref{u_delta_integral}) computed by MATHEMATICA
(see the MATHEMATICA notebook available at
 \url{https://shorturl.at/d0bgQ} accessed on 11 January 2026).
Also, the asymptotic representations for $\tau
\rightarrow 0$, i.e., (\ref{u_delta_tau->0}),\ and $\tau \rightarrow \infty $%
, i.e., (\ref{u_delta_tau->inf}) are numerically checked.

\begin{figure}[htbp]
 \includegraphics[width=0.8
\textwidth]{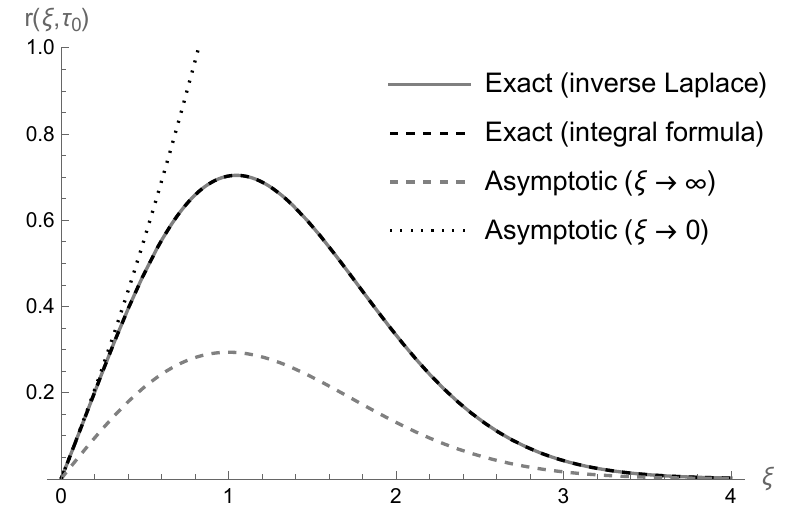}
\caption{Graph of $r\left( \protect\xi ,\protect\tau _{0}\right) $ with $%
\protect\tau _{0}=0.5$ for the delta pulse solutions and the asymptotic
approximations in the Kelvin--Voigt model. }
\label{Figure: u_delta_xi}
\end{figure}

Figure \ref{Figure: u_delta_xi} shows that the solution $r\left( \xi ,\tau
\right) $ given in terms of the inverse Laplace transform~(\ref{u_delta_inv_Laplace})\ is numerically equivalent to the integral solution (%
\ref{u_delta_integral}).
Indeed, the maximum discrepancy between both expressions in
Figure \ref{Figure: u_delta_xi} is $1.04126\times 10^{-11}$; thus, within the numerical error
of (\ref{u_delta_inv_Laplace}) and (\ref{u_delta_integral}) computed by MATHEMATICA
(see the MATHEMATICA notebook available at \url{https://shorturl.at/d0bgQ} accessed on 11 January 2026).
Also, the asymptotic representations for $\xi
\rightarrow \infty $, i.e., (\ref{u_delta_tau->0}), and $\xi \rightarrow 0$%
, i.e., (\ref{u_delta_xi->0}), are numerically checked.

\section{Conclusions}
Using the dimensionless variables $\xi$ and $\tau$, defined in (\ref{Xi_def}) and (\ref{Tau_def}), we have derived an integral-form expression for the mechanical response $r\left(\xi,\tau\right)$ of an initially quiescent, semi-infinite homogeneous Kelvin--Voigt medium to a pulse $r_{0}\left(\tau \right)$ applied at the origin. The general solution is given in (\ref{u(x,t)_general_dimensionless}),
\begin{equation*}
r\left( \xi ,\tau \right) =\int_{0}^{\tau }r_{0}\left( t_{\varepsilon
}\left( \tau -\tau ^{\prime }\right) \right) \,e^{-\tau ^{\prime }}\,g\left(
\xi ,\tau ^{\prime }\right) \,d\tau ^{\prime },
\end{equation*}
where the kernel $g\left( \xi ,\tau \right)$ is expressed in integral form by (\ref{g(xi,tau)_resultado})
\begin{equation*}
g\left( \xi ,\tau \right) =\frac{1}{4\sqrt{\pi }\,\tau ^{5/2}}\int_{\xi
}^{\infty }\left( v^{2}-2\tau \right) \exp \left( -\frac{v^{2}}{4\tau }%
\right) \,_{0}\overline{F}_{1}\left(
\begin{array}{c}
- \\
1%
\end{array}%
;\xi \left( v-\xi \right) \right) dv.
\end{equation*}%
For a delta-pulse excitation $r_{0}\left(\tau\right) = \delta\left( \tau \right)$, we obtained the explicit integral representation given in (\ref{u_delta_integral}):
\begin{equation*}
r\left( \xi ,\tau \right) =\frac{e^{-\tau }}{4\sqrt{\pi }\,\tau ^{5/2}}%
\int_{\xi }^{\infty }\left( v^{2}-2\tau \right) \exp \left( -\frac{v^{2}}{%
4\tau }\right) \,_{0}\overline{F}_{1}\left(
\begin{array}{c}
- \\
1%
\end{array}%
;\xi \left( v-\xi \right) \right) dv.
\end{equation*}%
Numerical comparisons confirm that this formulation is computationally efficient and fully equivalent to the corresponding inverse Laplace transform solution
given in (\ref{u_delta_inv_Laplace})
\begin{equation}
r\left( \xi ,\tau \right) =e^{-\tau }\,\mathcal{L}^{-1}\left[ \exp \left( -\,\xi \sqrt{s}+\frac{\xi
}{\sqrt{s}}\right) ;\tau \right] .  \notag
\end{equation}%
In contrast, the integral solution given by Hanin \cite{Hanin} is less efficient from a computational point of view
(see the MATHEMATICA notebook available at \url{https://shorturl.at/d0bgQ} accessed on 11 January 2026):%
\begin{eqnarray*}
r\left( \xi ,\tau \right)  &=&\frac{2\xi }{\pi \,\tau }\int_{0}^{1}e^{-2\tau
\,u^{2}}\cos \left( 2u\left[ \xi -\tau \sqrt{1-u^{2}}\right] \right) du \\
&&+\frac{e^{-2\tau }}{\pi }\int_{0}^{\infty }e^{-\tau \,u}\left[ \sin \left(
\frac{2+u}{\sqrt{1+u}}\xi \right) -\sin 2\xi \right] du.
\end{eqnarray*}%
The integral solution given by Dozio \cite{Dozio} in dimensionless variables reads as%
\begin{equation*}
r\left( \xi ,\tau \right) =e^{-\tau }\left[ \frac{\xi ^{2}}{2}%
+\frac{2}{\pi}\int_{0}^{\infty }u\,e^{-\tau \,u^{2}}\sin \left( \xi \left[ \frac{1}{u}+u%
\right] \right) du\right] .
\end{equation*}%
However, it seems that this solution is not correct from our numerical
experiments (see the MATHEMATICA notebook {available at}
 \url{https://shorturl.at/d0bgQ} accessed on 11 January 2026).
Nonetheless, the delta-pulse solution derived here allows us to obtain simple and explicit asymptotic formulas for the response as $\tau \rightarrow 0,\infty$ and $\xi \rightarrow 0,\infty $, summarized in (\ref{u_delta_tau->0})--(\ref{u_delta_xi->0}):
\begin{equation}
r\left( \xi ,\tau \right) \approx \frac{\xi }{2\sqrt{\pi }\tau ^{3/2}}\exp
\left( -\frac{\xi ^{2}}{4\tau }-\tau \right) ,\quad \tau \rightarrow 0 \quad \text{or} \quad
\xi \rightarrow \infty ,  \notag
\end{equation}%
\begin{equation}
r\left( \xi ,\tau \right) \approx \xi \,e^{-\tau }\left[ \frac{\xi }{2}+%
\frac{1}{\sqrt{\pi \,\tau }}\left( 1+\frac{1}{2\tau }\right) \right] ,\quad
\xi \rightarrow 0, \notag
\end{equation}%
and,
\begin{eqnarray*}
&&r\left( \xi ,\tau \right) \approx \xi \,e^{-\tau } \\
&&\times\left[ \frac{2\xi ^{2}\sqrt{\tau }}{9\sqrt{\pi }}\,_{0}F_{2}\left(
\begin{array}{c}
- \\
\frac{5}{2},\frac{5}{2}%
\end{array}%
;\frac{\xi ^{2}\tau }{4}\right) +\frac{\xi }{2}\,_{0}F_{2}\left(
\begin{array}{c}
- \\
\frac{3}{2},2%
\end{array}%
;\frac{\xi ^{2}\tau }{4}\right) +\frac{1}{\sqrt{\pi \,\tau }}%
\,_{0}F_{2}\left(
\begin{array}{c}
- \\
\frac{3}{2},\frac{3}{2}%
\end{array}%
;\frac{\xi ^{2}\tau }{4}\right) \right]    \\
&&\tau \rightarrow \infty  .
\end{eqnarray*}%

Similarly, we have obtained an integral-form solution for a step-pulse excitation $r_{0}\left(\tau\right) =\theta\left(\tau\right)$, given in (\ref{u(xi,tau)_theta_2})
\begin{eqnarray*}
r\left( \xi ,\tau \right)  &=&\int_{\xi }^{\infty }\,_{0}\overline{F}_{1}\left(
\begin{array}{c}
- \\
1%
\end{array}%
;\xi \left( v-\xi \right) \right) \\
&&\left[ \frac{\exp \left( -\frac{v^{2}}{4\tau }-\tau \right) }{\sqrt{\pi
\,\tau }}+\frac{e^{-v}}{2}\mathrm{erfc}\left( \frac{v-2\tau }{2\sqrt{\tau }}%
\right) -\frac{e^{v}}{2}\,\mathrm{erfc}\left( \frac{v+2\tau }{2\sqrt{\tau }}%
\right) \right] dv.
\end{eqnarray*}
This solution is numerically equivalent to the inverse Laplace transform representation \mbox{in (\ref{u(xi,tau)_Laplace})}
\begin{equation}
r\left( \xi ,\tau \right) =\mathcal{L} ^{-1}\left[ \frac{1}{s}\exp \left( -\,%
\frac{\xi \,s\,}{\sqrt{1\,+s}}\right) ;\tau \right] , \notag
\end{equation}%
and is more convenient for both numerical evaluation and asymptotic analysis than earlier formulations found in the literature.
For instance, the integral solution given by Morrison~\cite{Morrison},
\begin{equation*}
r\left( \xi ,\tau \right) =\frac{e^{-\tau}}{\pi}\int_{0}^{\tau }\frac{\cos \left( 2\sqrt{u\left(
\tau -u\right) }\right) }{\sqrt{u\left( \tau -u\right) }}\exp \left( 2u-%
\frac{\xi ^{2}}{4u}\right) du,
\end{equation*}
seems to be simpler, but it is numerically difficult to compute for $\tau \gg 1$ (see the MATHEMATICA notebook available at
 \url{https://shorturl.at/d0bgQ} accessed on 11 January 2026).
Also, it does not seem to be an easy task to obtain from Morrison's integral solution the asymptotic behavior as $\xi,\tau \rightarrow 0, \infty$.
Nevertheless, from our step-pulse integral representation, we derived asymptotic expressions for small and large values of $\xi$ and $\tau$, reported in
(\ref{u(xi,tau)_tau->0}), (\ref{u(xi,tau)_tau_inf2b}) and (\ref{u(xi,tau)_xi->0}):
\begin{equation}
r\left( \xi ,\tau \right) \approx \mathrm{erfc}\left( \frac{\xi }{2\sqrt{%
\tau }}\right) ,\quad \tau \rightarrow 0 \quad \text{or} \quad \xi \rightarrow \infty , \notag
\end{equation}
\begin{eqnarray*}
&&r\left( \xi ,\tau \right)
\approx 1-\xi ^{2}\,e^{-\tau }\left[ \frac{1}{2}\,_{0}F_{2}\left(
\begin{array}{c}
- \\
\frac{1}{2},2%
\end{array}%
;\frac{\xi ^{2}\tau }{4}\right) +\frac{2\xi \sqrt{\tau }}{3\sqrt{\pi }}%
\,_{0}F_{2}\left(
\begin{array}{c}
- \\
\frac{3}{2},\frac{5}{2}%
\end{array}%
;\frac{\xi ^{2}\tau }{4}\right) \right] ,   \\
&&\tau \rightarrow \infty ,
\end{eqnarray*}
and
\begin{equation}
r\left( \xi ,\tau \right) \approx 1-\frac{e^{-\tau }}{\sqrt{\pi \,\tau }}%
\,\xi -\frac{1+e^{-\tau }}{2}\,\xi ^{2},\quad \xi \rightarrow 0. \notag
\end{equation}

Finally, we note that all numerical checks and the graphical illustrations presented in this paper were carried out using MATHEMATICA. The corresponding MATHEMATICA notebook is available at
 \url{https://shorturl.at/d0bgQ} accessed on 11 January 2026.


\section*{Acknowledgments}

 The research activity of F. Mainardi and A. Mentrelli  has been carried out in the framework of the activities of the National Group of Mathematical Physics
 (GNFM, INdAM).
\\
 The authors are grateful to the anonymous referees for valuable suggestions which help us to improve the presentation of the results.

\section*{Funding}
A. Mentrelli is partially funded by the European Union – NextGenerationEU under the National Recovery and Resilience Plan (PNRR) - Mission 4 Education and research, Component 2 From research to business – Investment 1.1 Notice PRIN 2022 – DD N. 104 dated 2/2/2022, entitled ‘‘The Mathematics and Mechanics of Non-linear Wave Propagation in Solids’’ (proposal code: 2022P5R22 A; CUP: J53D23002350006), and by the Italian National Institute for Nuclear Physics (INFN), grant FLAG.

\section*{Author Contribution}
 Conceptualization: J.L. Gonz\'{a}lez-Santander and  F. Mainardi;
methodology: J.L. Gonz\'{a}lez-Santander,
 F. Mainardi, A. Mentrelli;
software: J.L. Gonz\'{a}lez-Santander;
data curation:   J.L. Gonz\'{a}lez-Santander,    F. Mainardi, A. Mentrelli;
writing---original draft preparation: J.L. Gonz\'{a}lez-Santander;
writing---review and editing:
J.L. Gonz\'{a}lez-Santander, F. Mainardi, A. Mentrelli;
supervision: F. Mainardi and  A. Mentrelli.
All authors have read and agreed to the published
version of the manuscript.

\section*{Conflicts of interests}
The authors declare no conflicts of interest.










\end{document}